\documentclass[a4paper,12pt]{article} 
\usepackage{amsmath, amsthm, amssymb}
\usepackage{url}

\usepackage{physics}
\usepackage{graphicx,lipsum}
\graphicspath{ {./downloads/} }

\usepackage[colorlinks,citecolor=red,urlcolor=blue,bookmarks=false,hypertexnames=true]{hyperref}
\usepackage{tikz}
\usetikzlibrary{calc}
\usetikzlibrary{shapes}
\usepackage[autostyle]{csquotes}
\makeatletter

\usepackage[colorlinks]{hyperref}
\usepackage[nameinlink,capitalize]{cleveref}
\newtheorem{theorem}{Theorem}[section]

\newtheorem{lemma}[theorem]{Lemma}
 \newtheorem{proposition}[theorem]{Proposition}
 \newtheorem{remark}[theorem]{Remark}
\newtheoremstyle{named}{}{}{\itshape}{}{\bfseries}{.}{.5em}{\thmnote{#3's }#1}
\theoremstyle{named}

\theoremstyle{definition}
\newtheorem{definition}[theorem]{Definition}
\newtheorem{example}[theorem]{Example}

\usepackage{xspace}
\usepackage[margin=1.0in]{geometry}

\usepackage{titlesec}

\usepackage{mathtools}

\DeclarePairedDelimiterX{\inp}[2]{\langle}{\rangle}{#1, #2}
\titleformat{\chapter}
  {\Large\bfseries} 
  {}                
  {0pt}            
  {\huge}

\begin{document}

\begin{center}
\fontsize{13pt}{10pt}\selectfont
    \textsc{\textbf{ON GRADED WEAKLY PRIMAL SUBMODULES OVER GRADED MODULES}}
    \end{center}
\vspace{0.1cm}
\begin{center}
   \fontsize{12pt}{10pt}\selectfont
    \textsc{{\footnotesize Tamem Al-shorma$n$, and  Malik Bataine$h$, }}

\end{center}
\vspace{0.2cm}

\begin{abstract}
Let $M$ be G-graded R-module. The idea of a graded weakly primal submodule of $M$, which is a generalization of a graded primal submodule, is introduced and discussed in this paper. Some characteristics and characterizations are assigned to these submodules and their homogeneous components.
\end{abstract}

\section{introduction}
Grading on a ring and its modules frequently facilitates computations by allowing one to concentrate on homogeneous elements, which are supposedly simpler or more manageable than random elements. However, in order for this to function, one must first understand how the structures are graded. One solution to this problem is to define the structures entirely in terms of G-graded R-modules, obviating the need to consider non-G-graded R-modules or non-homogeneous elements. 

All rings are assumed to be commutative with identity throughout this paper. The concept of weakly primal ideal was introduced by S. Atani and A. Darani in \cite{atani2007weakly}. In \cite{atani2009weakly}, S. Atani and A. Darani extended the concept of weakly primal ideals to modules. Let M be a R-module and N be a submodule of M. An element $x\in R$ is called weakly prime (simply wp) to N if $0\neq xm \in N$ for some $m \in M$, then $m \in N$. The submodule N is said to be a weakly primal submodule of M if the set $w(N)$ of elements of R that are not weakly prime to N forms an ideal of R, this ideal is always a weakly prime ideal of R, called the disjoint ideal P of N. We may also state that N is a P-weakly primal submodule in this situation. The concept of a graded weakly primal submodule of a G-graded R-module M over the G-graded ring R is defined in this paper, and subsequently, this class of submodules is studied.

We begin by defining the notations and terminologies that will be used throughout. Let G be an abelian group with identity e. A ring R is called a G-graded ring if $ R= \bigoplus\limits_{g \in G} R_g$   with the property $R_gR_h\subseteq R_{gh}$ for all $g,h \in G$, where $R_g$ is an additive subgroup of R for all $g\in G$. The elements of $R_g$ are called homogeneous of degree g. If $x\in R$, then $x$ can be written uniquely as $\sum\limits_{g\in G} x_g$, where $x_g$ is the component of $x$ in $R_g$. The set of all homogeneous elements of R is $h(R)= \bigcup\limits_{g\in G} R_g$. Let P be an ideal of a G-graded ring R. Then I is called a graded ideal if $I=\bigoplus\limits_{g\in G}I_g$, i.e, for $x\in I$ and  $x=\sum\limits_{g\in G} x_g$ where $x_g \in I_g$ for all $g\in G$. In \cite{atani2006grade}, I is called a graded weakly prime ideals of G-graded ring R if whenever $0\neq xy \in  I$, then we have that $x\in I$ or $y \in I$ for some $x,y \in h(R)$.

Let R be a G-graded ring . A left R-module M is said to be a graded R-module if there exists  a family of additive subgroups $\{M_g\}_{g \in G}$ of M such that $M= \bigoplus\limits_{g \in G} M_g$ with the property $M_gM_h\subseteq M_{gh}$ for all $g,h \in G$. The set of all homogeneous elements of M is $h(M)= \bigcup\limits_{g\in G} M_g$. Note that $M_g$ is an $R_e$-module for every $g\in G$. A submodule N of M is said to be graded submodule of M if $N=\bigoplus\limits_{g\in G}N_g$. Suppose that N is a graded submodule of M and I is a graded ideal of R. Then $(N:_RM)$ is defined as $(N:_RM) = \{r\in R: rM\subseteq N \}$ and then $(N:_RM)$ is a graded submodule of M. The graded submodule $(N:_MI)$ is defined as $(N:_MI)=\{m\in M : Im \subseteq N\}$. Particularly, we use $(N:_Ms)$ instead of $(N:_MRs)$.

The content of the paper is briefly summarized here. In Section 2, we give some basic results about graded weakly primal submodules. Every graded weakly prime submodule and every graded weakly primary submodule is a graded weakly primal submodule in Theorem \ref{thm4}. In Proposition \ref{prop1}, we give a description of graded weakly primal submodules. It is shown in Theorem \ref{thm2} that is N is a  graded weakly primal submodule of M, then GW(N) is a graded weakly prime ideal of R. In Theorem \ref{thm5} we prove that if M is a faithful G-graded finitely generated multiplication R-module and if R is a G-graded WP-ring, then M is a G-graded WP-module. The concepts graded primal submodule and graded weakly
primal submodule are different. In fact, neither implies the other (see Examples \ref{exm1}). In Theorem \ref{thm3} we prove that if N be a graded weakly primal submodule of M with $(N:_RM) \subseteq P$ and $N(N:_RM) \neq 0$, then N is a graded primal submodule of M.

Let $R$ is a G-graded ring and that $S$ is a multiplicative closed subset of $h(R)$. Consider the G-graded $R_S$-module $M_S$. In Section 3, the weakly primal submodules of M and the weakly primal submodules of $M_S$ are studied with respect to each other. It is shown in Theorem \ref{thm8} that there is  a one-to-one correspondence between the graded P-weakly primal submodule of M and the graded $P_S$-weakly primal submodule of $M_S$ in which P is a graded weakly prime ideal of R, and $S\subseteq h(R)$ is a multiplucative closed set of R with $P\cap S =\emptyset$.

\section{Graded weakly primal submodules}
To begin, consider the following definitions:
\begin{definition}\label{def1}
Let R be a G-graded ring and P be a graded ideals of R.
\\

(1) The element $x \in h(R)$ is said to be a graded weakly prime (simple $gwp$) to P if $0\neq xy \in N$ for some $y \in h(R)$, then $y \in P$.
\\

(2) The element $x \in h(R)$ is not graded weakly prime (simply $ngwp$) to N if $0\neq xy \in P$ for some $y\in h(R)$, then $y \in R/P$.
\\

(3) The set of all elements of h(R) that are not graded weakly prime ($ngwp$) to P is denoted by gw(N).
\\

(4) The graded ideal P of R is called a graded weakly primal ideal of R if gw(N) is a graded ideal of R.
\end{definition}

\begin{definition}\label{def2}
Let M be a G-graded R-module and N be a graded submodule of M.
\\

(1) The element $x \in h(R)$ is said to be a graded weakly prime (simple $GWP$) to N if $0\neq xm \in N$ for some $m \in h(M)$, then $m \in N$.
\\

(2) The element $x \in h(R)$ is not graded weakly prime (simply $NGWP$) to N if $0\neq xm \in N$ for some $m\in h(M)$, then $m \in M/N$.
\\

(3) The set of all elements of h(R) that are not graded weakly prime ($NGWP$) to N is denoted by GW(N).
\\

(4) The graded submodule N of M is called a graded weakly primal submodule of M if GW(N) is a graded ideal of R.
\end{definition}

\begin{lemma}\label{lem1}
Let M be a G-graded R-module and N be a graded submodule of M. If N is a graded primal submodule of M, then N is a graded weakly primal submodule of M.
\\
\\
\textbf{Proof.} Let N be a graded submodule of M that is graded primal. Then $G(N)$ is a graded ideal of R. We must now prove that N is a graded weakly primal submodule of M, and that $GW(N)$ is a graded ideal of R. Since $GW(N)$ is a graded ideal of R, N is a graded weakly primal submodule of M if $GW(P) = R$. Now, if $GW(N) \neq R$, then $1\notin GW(N)$ is true. $GW(N) = G(N)$, we show. It is clearly, $G(N) \subseteq GW(N)$. Conversely, let $w \in GW(N)$ then there exists $x \in h(R)$ with $ 0\neq xw\in N$ as a result $x \notin N$. If $x \in N$, then $1 \in GW(N)$ and then $GW(N) = R$ which is a contradiction. So, $x \notin N$ and hence $w \in G(N)$. Thus $G(N) \subseteq GW(N)$ then $G(N) = GW(N)$. Therefore, N is a graded weakly primal submodule of M.
\end{lemma}

\begin{remark}\label{rem1}
Let M be a G-graded R-module and let N be a graded submodule of M. Then the following statements are true:
\\

(1) The element $0 \in h(R)$ is invariably $GWP$ to N.
\\

(2) If $x \in h(R)$ is a graded prime to N, then x is a $GWP$ to N.
\\

(3) The set of all not weakly prime to N (W(N)) is subset of the set of all not graded weakly prime to N (GW(N)).

\end{remark}

If $N$ is a graded primal submodule of $M$, then $N$ is a graded weakly primal submodule of $M$, as proven in the Lemma \ref{lem1}.  However, the converse is not always true. In addition, the converse is not true in (2) of the Remark \ref{rem1}.

\begin{example}\label{exm1}
(1) Let $R= \mathbf{Z}$ be a G-graded ring with $R_0 = \mathbf{Z}$ and $R_g = \{0\}$ for all $g\in G-\{e\}$. Take R as a G-graded R-module. Then $N=12\mathbf{Z}$ is a graded submodule of R. Then N is a graded weakly primal submodule of R with $GW(N) = \mathbf{Z}$, since $0\neq 3.4=12 \in N$ and $3 \not\in N$, $4 \notin N$. On the other hand, N is not a graded primal submodule of R.
\\

(2) Let $M=\mathbf{Z}/ 24 \mathbf{Z}$ be a G-graded $\mathbf{Z}$-module where $M_0 = M$ and $M_g = \{0\}$ for all $g\in G$ and $N = 8\mathbf{Z}/24\mathbf{Z}$ be a graded submodule of M. Since N is a primal submodule by \cite[Example 2.10]{atani2009weakly}, then N is a graded primal submodule of M. Now, since $\Bar{2}, \Bar{4} \in M/N$ with $0\neq 2.\Bar{4} \in N$ and $0\neq 4. \Bar{2} \in N$, then $2,4 \in Gw(N)$. If $6. \Bar{x} \in N$ for some $\Bar{x} \in M$, then $4$ divides x and then $6\Bar{x}=0$. This prove that $2+4=6$ is a GWP to N. Thus $GW(N)$ is not graded ideal of $R=\mathbf{Z}$. Therefore, N is not graded weakly primal submodule of M.  
\\

(3) Let $M = \mathbf{Z}_{12}$ be a G-graded $\mathbf{Z}$-module and $N = \{\Bar{0}\}$ be a graded submodule of M. It is clearly $GW(N) = \emptyset$, then N is a graded weakly primal submodule of M. Now, since $3.\Bar{4} = \Bar{0} \in N$ and $4.\Bar{3}=\Bar{0}\in N$, then we have that $3,4 \in G(N)$, while $4-3=1$ is a graded prime to N. Therefore, N is not graded primal submodule of M.
\\

(4) Let $M=\mathbf{Z}/ 32\mathbf{Z}$ be a G-graded $\mathbf{Z}$-module where $M_0 = M$ and $M_g = \{0\}$ for all $g\in G$ and $N = 8\mathbf{Z}/32\mathbf{Z}$ be a graded submodule of M. Then $4.8=0\in N$, 4 is not graded prime to N. But if $4. \Bar{x} \in N$ where $\Bar{x}$ is the coset in M, then 4 is divides of x. Thus $4.\Bar{x} = \Bar{0}$. Therefore 4 is $GWP$ to N.
\end{example}

\begin{lemma}\label{lem2}
Let M be a G-graded R-module and N be a graded submodule of M, then:
\\

(1) $(N:_RM)$ is subset of the set of all $NGWP$ to N (GW(N)).
\\

(2) The set of all $ngwp$ to $(N:_RM)$ $(gw((N:_RM)))$ is subset of the set of all $NGWP$ to N (GW(N)).
\\
\\
\textbf{Proof.} (1) Assume that $x \in h(R)$ such that $x \in (N:_RM)$. $(N:_RM)$ is a graded ideal of R since N is a graded submodule of M. Take $x \in (N:_RM)$ we have that x is not weakly prime to N. Thus by (3) in Remark \ref{rem1}, $x \in W(N)\subseteq GW(N)$. Therefore, $(N:_RM) \subseteq GW(N)$.
\\

(2) Let $x \in h(R) \cap gw((N:_RM))$, so there exists  $y \in h(R)-(N:_RM)$ with $0\neq xy \in (N:_RM)$, there exist $m \in h(M)$ with $ym \notin N$. It therefore follows that $xym \in N$ with $ym \notin N$ we have that x is not weakly prime to N. Thus $x \in W(N)\subseteq GW(N)$ by (3) in Remark \ref{rem1}. Therefore, $gw((N:_RM)) \subseteq GW(N)$.
\end{lemma}

\begin{definition}\label{def3}
Let M be a G-graded R-module, N be a graded weakly primal submodule of M and P be a graded ideal of R. Then P is called the graded weakly adjoint ideal of N, and also we said to be N is a graded P-weakly primal submodule of M.  
\end{definition}

A characterization of graded weakly primal submodules is provided by the following Proposition.

\begin{proposition}\label{prop1}
Let M be a G-graded R-module, N be a graded submodule of M and P be a graded ideal of R. Then the following statements are equivalent:
\\

(1) N is a graded P-weakly primal submodule of M.
\\

(2) For every $p \notin P-\{0\}$, $(N:_Mp)= N\cup (0:_Mp)$, and for every $0\neq p \in P$, $N \cup (0:_Mr) \subset (N:_Mp)$.
\\
\\
\textbf{Proof.} (1) $\Rightarrow$ (2) Assume that N is a graded P-weakly primal submodule of M. Let $p \notin GW(N)=P-\{0\}$ and $m \in (N:_Mp)$. If $pm =0$, then $m \in (0:_Mp)$. Now, suppose that $pm \neq 0$, since p is a GWP to N we have that $m \in N$. Thus $m \in N\cup (0:_Mp)$, that is $(N:_Mp)\subseteq N\cup (0:_Mp)$. Therefore, $(N:_Mp) =N\cup (0:_Mp)$. Now, suppose that $p \in GW(N)$. Thus p is NGWP to N, then there exist $m \in h(M)$ such that $0\neq pm \in N$ with $m \in M/N$. Hence $m \in (N:_Mp)-(N\cup (0:_Mp))$.
\\

(2) $\Rightarrow$ (1) $GW(N) = P-\{0\}$ follows from (2). As a result, N is a graded P-weakly primal submodule of M.
\end{proposition}

\begin{theorem}\label{thm1}
Let M be a cyclic G-graded R-module and N be a graded submodule of M. If N is a graded weakly primal submodule of M, then $(N:_RM)$ is a graded weakly primal ideal of R.
\\
\\
\textbf{Proof.} Suppose that $M = Rm$ for any $m \in h(M)$, and we set $P=(N:_RM)$. We prove that $gw(P)=GW(N)$. Let $x \in gw(P)$, so there exists $r \in R/P$ with $0\neq xr \in P$. In this instance, $xrm \neq 0$ because otherwise $xr \in (0:_Rx)=(0:_RM)=0$, a contradiction. As $rm \in M/N$ we have that x is $NGWP$ to N, thus $x \in GW(N)$ and thus $gw(P) \subseteq GW(N)$. Now, suppose that $x \in GW(N)$, so $0\neq xm' \in N$ for some $m' \in M/N$. But we can write $m' = x'm$ for some $x' \in h(R)$. Thus $0\neq xx'm\in N$, this implies that $0\neq xx'\in P$ with $x'\in R/P$, then x is ngwp to P, then $x \in gw(P)$ and then $GW(N) \subseteq gw(P)$. Therefore, $gw(P)=GW(N)$. 
\end{theorem}

\begin{theorem}\label{thm2}
Let M be a  G-graded R-module and N be a graded submodule of M. If N is a graded weakly primal submodule of M, then GW(N) is a graded weakly prime ideal of R.
\\
\\
\textbf{Proof.} Assume that $0\neq xy \in GW(N)$ for some $x,y \in h(R)$ but $x\notin GW(N)$. Then $xy$ is NGWP to N  and then $0\neq xym \in N$ for some $m \in h(M)$. As $x\notin N $, so x is a GWP to N. Thus we have that $x(ym)\in N$, so $ym \in N$. Hence y is nwp to N. Thus $y \in W(N) \subseteq GW(N)$ by Remark \ref{rem1}. Therefore, GW(N) is a graded weakly prime ideal of R. 
\end{theorem}

\begin{theorem}\label{thm3}
Let M be a G-graded R-module and  N be a graded weakly primal submodule of M. If $(N:_RM) \subseteq P$ and $N(N:_RM) \neq 0$, then N is a graded primal submodule of M.
\\
\\
\textbf{Proof.} It is enough to show that GW(N) = G(N). Let $0\neq x \in GW(N)$, x is NGWP to N. Thus by Remark \ref{rem1} (2), x is not graded prime to N, then we have $x \in G(N)$. Thus $GW(N) \subseteq G(N)$. Now, we need to show that the converse $G(N) \subseteq GW(N)$, let $y \in G(N)$ so for some $m\in h(M)$ then $m \in M/N$ with $ym \in N$. If $ym \neq 0$, then y is  NGWP to N and then $y \in GW(N)$. But if $ym = 0$ we have tow cases:
\\

\textbf{Case one:} $yN \neq 0$. Then there exists $z \in N$ with $yz \neq 0$. Now $0\neq y (m+z) \in N$ with $m+n \in M/N$, then we have that y is NGWP to N and then $y \in GW(N)$.
\\

\textbf{Case two:} $yN =0$. If $m(N:_RM) \neq 0$, then $zm \neq 0$ for some $z \in (N:_RM)$. Now, $0\neq (y+z)m \in N$ with $m \in M/N$ we have that $y+z$ is NGWP to N, then $y+z \in GW(N)$ and then $y \in GW(N)$. If $m(N:_RM)=0$, but $N(N:_RM) \neq 0$, then there is $x \in (N:_RM)$ and $n \in N$ with $xn \neq 0$. Now, $0 \neq (y+x)(m+n) \in N$ with $m+n \in M/N$, then $y+x$ is NGWP to N. Thus $y+x \in GW(N)$ with $x \in GW(N)$. Therefore $G(N) \subseteq GW(N)$.
\\
Thus $GW(N)=G(N)$ which is implies that N is a graded primal submodule of M.
\end{theorem}

Recall that a graded submodule N of A G-graded R-module M is called a graded weakly prime submodule of M if $N \neq M$, and whenever $x \in h(R)$ and $m \in h(M)$ with $0\neq xm \in N$, then either $m \in N$ or $x \in (N:_RM)$ (see \cite{atani2006graded}).

\begin{definition}\label{def4}
Let M be a G-graded R-module and N be a graded submodule of M. N is called a graded weakly primary submodule of M if whenever $x\in h(R)$ and $m \in h(M)$ with $0\neq xm \in N$, then either $m \in N $ or $x^n \in (N:_RM) $ for some positive integer n.
\end{definition}

\begin{theorem}\label{thm4}
Let M be a G-graded R-module and N be a graded submodule of M. Then the following statements are hold:
\\

(1) If N is a graded weakly primary submodule of M, then N is a graded weakly primal submodule of M. 
\\

(2) If N is a graded weakly prime submodule of M, then N is a graded weakly primal submodule of M.
\\
\\
\textbf{Proof.} (1) Let N be a graded primary submodule of M. we claim that $Grad((N:_R:M))=GW(N)$. Let $x \in Grad((N:_RM))$. If $x \in (N:_RM)$, $x \in GW(N)$ by (1) in Lemma \ref{lem2}. Now, suppose that $x \in Grad(N)-(N:_RM)$. If $x \in (N:_RM)$, then $x \in W(N)$. If $x \notin (N:_RM)$, then there exist a positive integer $n > 1$ for which $x^n \in (N:_RM)$. In this case $xx^{n-1} \in (N:_RM)$ with $x^{n-1} \in h(R)-(N:_RM)$ which is implies that $x \in W((N:_RM))$. Thus $x \in gw((N:_RM)) \subseteq GW(N)$ by (2) in Lemma \ref{lem2}, that is $Grad((N:_RM)) \subseteq GW(N)$. Conversely, suppose that $y \in GW(N)$ and there exist $m \in h(M)-N$ with $y \in N$. As N is a graded weakly primary submodule of M we have that $y \in Grad((N:_RM))$. Thus $GW(N) \subseteq Grad((N:_RM))$. Therefore, $Grad((N:_RM))=GW(N)$ and thus N is a graded weakly primal submodule of M.
\\

(2) Because any graded weakly prime submodule is also a graded weakly primary submodule, (1) follows.
\end{theorem}

Recall that, a $G$-graded $R$-module $M$ is called a $G$-graded multiplication $R$-module if $N=MI$ for some graded ideal $I$ of $R$ for every graded submodule $N$ of $M$ (see \cite{escoriza1998multiplication}). If $N$ is a graded submodule of a $G$-graded multiplication $R$-module $M$, then $N = (N:_R M)M$ is simple to prove. A $G$-graded $R$-module $M$ is said to be graded finitely generated if $M = \sum\limits_{i=1}^{n} Rm_i$, where $m_i \in h(M)$.

\begin{remark}\label{rem2} \cite[Remark 1.7]{darani2011graded}
If M is a G-graded finitely generated multiplication R-module and P is a graded ideal of R containing $(0 :_R M)$, then $(PM :_R M) =P$.
\end{remark}

\begin{proposition}\label{prop2}
Let M be a G-graded finitely generated multiplication R-module and P be a graded ideal of R. If P is a graded weakly primal ideal of R containing $(0:_RM)$, then PM is a graded weakly primal ideal of M.
\\
\\
\textbf{Proof.} It is clearly $gw(P)=gw((PM:_RM)) \subseteq GW(PM)$ by (2) Lemma \ref{lem2}. Now, suppose that $x \in h(R)$ is NGWP to PM, then there exist $m \in h(M)-PM$ with $0\neq xm \in PM$. Then we have $RxRm \subseteq PM$. Since M is a G-graded multiplication R-module, then there exist a graded ideal I of R with $Rm=IM$. Then we have that $(RxI)M \subseteq PM$ and then $RxI \subseteq (PM:_RM)=P$ by Remark \ref{rem2}. If $I \subseteq P$, then $m \in Rm = IM \subseteq PM$ which is a contradiction. Then there exists $y \in I-P$ where $y \in h(R)$. In this case $xy \in RxJ \subseteq P$ which is implies that x is ngwp to P. Therefore is $x \in gw(P)$. Thus $GW(PM) \subseteq gw(P)=gw(PM:_RM)$. Therefore, PM is a graded weakly primal submodule of M.
\end{proposition}

\begin{definition}\label{def5}
Let M be a G-graded ring, and M be a G-graded R-module.
\\

(1) A G-graded ring R is said to be a G-graded WP-ring if every graded ideal of R is a finite product of graded weakly primal ideal of R.
\\

(2) A G-graded R-module M is said to be a G-graded WP-module if every graded submodule N of M has a graded weakly primal factorization $N=P_1P_2...P_nN^\ast$ where $P-1,P_2,...,P_n$ is a graded weakly primal ideals of R and $N^\ast$ is a graded weakly primal submodule of M.
\end{definition}

\begin{theorem}\label{thm5}
Let R be a G-graded ring and M be a faithful G-graded finitely generated multiplication R-module. If R is a G-graded WP-ring, then M is a G-graded WP-module.
\\
\\
\textbf{Proof.} Suppose that R is a G-graded WP-ring and N is a graded submodule of M. Then N=PM for some graded ideal P of R. Since R is a G-graded WP-ring, P has a factorization $P=P_1P_2...P_n$ where $P_i$ is a graded weakly primal ideal of R for all $i=\{1,2,..n\}$. In this case $N=PM=(P_1P_2...P_{n-1})(P_nM)$. By Proposition \ref{prop2}, $P_nM$ is a graded weakly primal submodule of M. Therefore, M is a G-graded WP-module.
\end{theorem}

\begin{proposition}\label{prop3}
Let M be a faithful G-graded R-module and N be a graded submodule of M. Then the following holds:
\\

(1) If N is a graded P-weakly primal submodule of M, then $(N:_RM) \subseteq P$.
\\

(2) If N is a graded 0-weakly primal submodule of M, then M/N is a faithful G-graded R-module.
\\
\\
\textbf{Proof.} (1) Suppose that N is a graded P-weakly primal submodule of M. Since N is a graded submodule of M, so $(N:_RM)$ is a graded ideal of R. There exists $m \in h(M/N)$ and $x \in (N:_RM)$, then $xm \in N$. If $xm \neq 0$, then x is NGWP to N, that $x \in P$. If $xm =0$, then there exists $n \in N$ with $xn \neq 0$. Thus $0\neq x(m+n)\in N$ with $(m+n)\notin N$ which is implies that x is NGWP to N, that is $x \in P$. Therefore, $(N:_RM) \subseteq P$.
\\

(2) As a result of (1)
\end{proposition}

\section{Graded modules of fractions}
This section will discuss the relations between the classes of graded weakly primal submodules of M and graded weakly primal submodules of $M_S$.

Assume that $R$ is a G-graded ring and that $S$ is a multiplicative closed subset of $h(R)$. Remember that if $R$ is a G-graded ring, then $R_S = \{ \frac{r}{s}: r \in h(R), s \in S \}$ is also a G-graded ring. In addition, if $P$ is a graded ideal of R, $P_S$ is a graded ideal of $R_S$. Now, let M be a G-graded R-module and consider the module fractions $M_S= \{ \frac{m}{s}: m \in h(M), s \in S \}$. Then $M_S$ is a G-graded $R_S$-module. 

\begin{theorem}\label{thm6}
Let M be a G-graded R-module and $S\subseteq h(R)$ be a multiplicatively closed set of R. Let N be a graded weakly primal submodule of M with $GW(N) \cap S =\emptyset$. Then the following holds:
\\

(1) If $0\neq \frac{n}{s} \in N_S$, then $n \in N$
\\

(2) If L is a graded submodule of M, then $(N:_RL)_S = (N_S:_{R_S}L_S)$.
\\
\\
\textbf{Proof.} (1) Since $\frac{n}{s}\in N_S$, so $\frac{n}{s}\in h(M_S)$. Now, there exists $x \in N$ and $t \in S$ such that $\frac{n}{s}=\frac{x}{t}$. Thus $0\neq wtn=wsx \in N$ for some $w \in S$. If $n \notin N$, then wt is NGWP to N and then we have $ut \in GW(N) \cap S$, which is a contradiction. Therefore, $n \in N$.
\\

(2) It is clearly $(N:_RL)_S \subseteq (N_S:_{R_S}L_S)$. Conversely, suppose that $\frac{x}{s}\in (N_S:_{R_S}L_S)$. So, for any $m \in L$, $\frac{xm}{s}=\frac{x}{s}.\frac{m}{1} \in N_S$. If $\frac{xm}{s}=0$, then there exists $t \in S$ with $txm =0$. Thus $xm =0 \in N$. If $\frac{xm}{0}\neq 0$, then $xm \in N$ by (1). Hence $x \in (N:_RL)$, so $\frac{x}{s} \in (N:_RL)_S$. Thus $(N_S:_{R_S}L_S) \subseteq (N:_RL)_S$. Therefore, $(N:_RL)_S = (N_S:_{R_S}L_S)$.
\end{theorem}

Let M be a G-graded R-module and $S\subseteq h(R)$ be a multiplucative closed set of R. Consider the homomorphism $\phi:M \rightarrow M_S$ which is defined by $\phi(m)=\frac{m}{1}$. Then $\phi$ is a homogeneous homomorphism of degree e. If N is a graded submodule of $M_S$, then we defined $N\cap M =\phi^{-1}(N)$.

\begin{proposition}\label{prop4}
Let M be a G-graded R-module and $S\subseteq h(R)$ be a multiplicative closed set of R. If N is a graded P-weakly primal submodule of $M_S$, then $N\cap M$ is a graded $P\cap R$-weakly primal submodule of M.
\\
\\
\textbf{Proof.} Since $P\cap R$ is a graded prime ideal of R, then $P\cap R$ is a graded weakly prime ideal of R. Now, it is enough to prove that $GW(N\cap M)=(P\cap R)-\{0\}$. Let $x \in GW(N\cap M)$, so there exist $m \in M-(N\cap M)$ such then $0\neq xm \in N\cap M$. As a result of $0\neq \frac{x}{1}.\frac{m}{1}\in N$ and $\frac{m}{1}\in M_S-N$ that $\frac{x}{1}\in P-\{0\}$. Hence $x \in P\cap R$ and then $GW(N\cap M) \subseteq (P\cap R)-\{0\}$. For the other direction, let $y \in (P\cap R)-\{0\}$, then we have that $0\neq \frac{y}{1}\in P$. Then there exists $\frac{m}{s}\in M_S-N$ with $0\neq \frac{y}{1}.\frac{m}{s}\in N$. Thus $0\neq \frac{ym}{1}\in N$. Then $0\neq ym \in N\cap M$ with $m \in M-(N\cap M)$ and then $y \in GW(N\cap M)$, hence $(P\cap R)-\{0\} \subseteq GW(N\cap M)$. Therefore, $N\cap M$ is a graded $P\cap R$-weakly primal submodule of M.
\end{proposition}

\begin{theorem}\label{thm7}
Let M be a G-graded R-module, N be a graded P-weakly primal submodule of M, and $S\subseteq h(R)$ be a multiplucative closed set of R with $P\cap S = \emptyset$. Then the following holds:
\\

(1) $N_S$ is a graded $P_S$-weakly primal submodule of $M_S$.
\\

(2) $N=(N_S\cap M)$.
\\
\\
\textbf{Proof.} (1) Let $\frac{0}{1}\neq \frac{x}{s}\in P_S$, then $0\neq x \in P$. Hence there exists $m \in h(M)$ such that $0\neq xm \in N$. If $0\neq \frac{x}{s}.\frac{m}{1}\in N_S$, and by Theorem \ref{thm6} , $\frac{m}{1}\notin N_S$. Thus $\frac{x}{s}$ is a NGWP to $N_S$, then we have $\frac{x}{a}\in GW(N_S)$. Then $P_S-\{0\} \subseteq GW(N_S)$. Now, assume that $\frac{y}{s}\in GW(N_S)$. So $0\neq \frac{y}{s}.\frac{m}{t}\in N_S$ for some $\frac{m}{t}\in h(M_S)$. Hence, by Theorem \ref{thm6} we have that $0\neq ym \in N$ with $m \in h(M)$ and then $y \in P-\{0\}$. Thus $\frac{y}{s}\in P_S-\{0\}$, then $GW(N_S)\subseteq P_S-\{0\}$. Then $GW(N_S)\cup \{0\} = P_S$. Therefore, $N_S$ is a graded $P_S$-weakly primal submodule of $M_S$.
\\

(2) It is clearly $N \subseteq (N_S\cap M)$. Suppose that $m \in N_S\cap M$. If $m=0$, then $m \in N$. If $M\neq 0$, then $0\neq \frac{m}{1} \in N_S$ we have that $m \in N$ bt Theorem \ref{thm6}. Therefore, $(N_S\cap M) \subseteq N$.
\end{theorem}

\begin{theorem}\label{thm8}
Let M be a G-graded R-module, P be a graded weakly prime ideal of R, and $S\subseteq h(R)$ be a multiplucative closed set of R with $P\cap S =\emptyset$. Then there exist a one-to-one correspondence between the graded P-weakly primal submodule of M and the graded $P_S$-weakly primal submodule of $M_S$.
\\
\\
\textbf{Proof.} This follows from Proposition \ref{prop4} and Theorem \ref{thm7}.
\end{theorem}

\section{Conclusions}
In this study, we introduced the concept of graded weakly primal submodules which is a generalization of graded weakly primal submodules. We investigated some basic properties of graded weakly primal submodules. As a proposal to further the work on the topic, we are going to study the concepts of graded S-primal submodules and graded weakly S-primal submodules. Also, we will also generalize the primal and Weakly Primal SubSemiModules (see \cite{bataineh2014primal}) on G-graded R-module and generalizations of primal ideals over commutative semirings (see \cite{bataineh2014}) on G-graded rings.

\bibliographystyle{amsplain}

\end{document}